\documentclass[12pt]{amsart}
\usepackage{amssymb,amsmath,amsthm,amscd}
\usepackage[all]{xy}

\numberwithin{equation}{section}

\newtheoremstyle{fancy1}{10pt}{10pt}{\itshape}{12pt}{\textsc\bgroup}{.\egroup}{8pt}{
}
\newtheoremstyle{fancy2}{10pt}{10pt}{}{12pt}{\itshape}{.}{8pt}{ }

\theoremstyle{fancy1}

\newtheorem{cor}[equation]{Corollary}
\newtheorem*{cora*}{Corollary A}
\newtheorem*{corb*}{Corollary B}
\newtheorem*{corc*}{Corollary C}
\newtheorem*{cord*}{Corollary D}
\newtheorem*{core*}{Corollary E}
\newtheorem{lem}[equation]{Lemma}
\newtheorem{prop}[equation]{Proposition}
\newtheorem{thm}[equation]{Theorem}

\newtheorem*{main*}{Result}
\newtheorem*{conjecture*}{Conjecture}
\newtheorem*{cor*}{Corollary}

\newcommand{\pa}[1]{\pi_A(U_{#1})}
\newcommand{\gal}[2]{\operatorname{Gal}\left( #1 / #2 \right)}

\newcommand{\dis}[2]{|\mathfrak{d}_{#1/\Q}|^{1/#2}}

\newtheorem*{def*}{Definition}

\newtheorem*{rem*}{Remark}
\newtheorem{remarks}[equation]{Remarks}

\newtheorem*{example*}{Example}
\newtheorem{examples}[equation]{Examples}
\newtheorem*{examples*}{Examples}

\newtheorem*{prop*}{Proposition}
\theoremstyle{remark}
\newtheorem*{case*}{Case}
\newtheorem*{case}{Case}
\newtheorem*{proofA*}{Proof of Corollary A}
\newtheorem*{proofB*}{Proof of Corollary B}
\newtheorem*{proofC*}{Proof of Corollary C}
\newtheorem*{proofD*}{Proof of Corollary D}
\newtheorem*{proofE*}{Proof of Corollary E}
\newtheorem*{proofprop*}{Proof of Proposition}
\newtheorem*{proofthm*}{Proof of Theorem \ref{maintheorem}}

\newcommand{\cref}[1]{Corollary~\ref{#1}}

\newcommand{\F}{{\mathbb{F}}}
\newcommand{\PP}{{\mathbb{P}}}
\newcommand{\A}{{\mathbb{A}}}
\newcommand{\Q}{{\mathbb{Q}}}

\newcommand{\Z}{{\mathbb{Z}}}
\def\con#1=#2(#3){#1 \equiv #2 \bmod{#3}}

\newcommand{\Aut}{\ensuremath{\operatorname{Aut}}}

\newcommand{\SL}{\ensuremath{\operatorname{SL}}}

\newcommand{\Spec}{\ensuremath{\operatorname{Spec}}}

\begin{document}

\title[Galois extensions ramified only at one prime]{Galois extensions ramified only at one prime}

\author{Jing Long Hoelscher}
\address{University of Arizona\\
     Tucson, AZ 85721}
\email{jlong@math.arizona.edu}

\begin{abstract}
This paper gives some restrictions on finite groups that can occur as Galois groups of extensions over $\Q$ and over $\F_q(t)$ ramified only at one finite prime. Over $\Q$, we strengthen results of Jensen and Yui about dihedral extensions and rule out some non-solvable groups. Over $\F_q(t)$ restrictions are given for symmetric groups and dihedral groups to occur as tamely ramified extension over $\F_q(t)$ ramified only at one prime.
\end{abstract}

\maketitle

%-------------- Article Text--------------------

\section*{Introduction}
This paper studies Galois groups with prescribed ramification in both the function field and number field cases. We are particularly interested in the case with a single finite ramified place. 

In the geometric case, we are concerned with curves over a field $k$ of characteristic $p>0$. Let $X$ be a smooth connected projective curve of genus $g$ over $k$ and let $S=\{\xi_1,...,\xi_n\}$ be a finite set of $n>0$ closed points on $X$. Then $U=X-S$ is an open subset of $X$. Define $\pi_A(U)$ to be the set of finite groups that occur as Galois groups over $X$ with ramifications only at $S$, and $\pi_A^t(U)$ the subset of $\pi_A(U)$ corresponding to covers in which only tame ramifications occur. In the case $k$ is algebraically closed, Corollary $2.12$ of Chapter XIII in \cite{Gro} implies that: if $G$ is a Galois cover of $X$ ramified only at $S$, then $G/p(G)$ has $2g+n-1$ generators. Here $p(G)$ denotes the quasi $p$-part of $G$, i.e. the subgroup of $G$ generated by elements of order a power of $p$. This statement can be carried over to the case where $k$ is a finite field $\F_q$ of order $q$, a power of $p$, if we restrict to regular covers $Y/X$, i.e. where $k$ is algebraically closed in the function field of $Y$. If $X/\F_q(t)$ is a tame regular cover with Galois group $G$, then $G$ has at most $2g+n-1$ generators. Here we count the number $n$ of ramified primes according to their degree, i.e. $n$ is the degree of $S$ as a divisor over $\mathbb{F}_q$.
Proposition \ref{conjugaterelation}, of Section $3$ below, gives further restrictions on $\pi_A^t(\mathbb{A}^1_{\F_q}-(f))$, where $(f)$ is the divisor of zeroes of an irreducible $f\in\mathbb{F}_q[t]$. The following two Corollaries of Proposition \ref{conjugaterelation} show that dihedral groups and symmetric groups tend to not occur in $\pi_A^t(\mathbb{P}^1_{\F_q}-(f))$.
\begin{cora*}For any integer $k\geq 1$ and irreducible $f\in\F_q[t]$, the dihedral group $D_{4k}\notin \pi_A^t(\mathbb{P}^1_{\F_q}-(f))$. If the degree $d=deg(f)$ is odd, we also have $D_{4k+2}\notin \pi_A^t(\mathbb{P}^1_{\F_q}-(f))$.
\end{cora*}
\begin{corb*}If $2|q$ and $n>2$, the symmetric group $S_n\notin\pi_A^t(\mathbb{P}^1_{\F_q}-(f))$ for each irreducible $f\in\F_q[t]$. If $2\nmid q$ and the prime $f$ is of odd degree, then the same conclusion holds.
\end{corb*}

In the arithmetic case, we consider Galois extensions over $\Q$. Denote $U_n=\Spec(\Z[\frac{1}{n}])$, an open subset of $\Spec(\Z)$. Then $\pi_A(U_n)$ is the set of finite groups that occur as Galois groups over $\Q$ ramified only at primes dividing $n$. Motivated by Corollary $2.12$ in Chapter XIII of \cite{Gro} and the analogy between function fields and number fields, Harbater posed a corresponding conjecture in \cite{Ha}:
\begin{conjecture*}\label{DavidConjecture}[Harbater, 1994]
There is a constant $C$ such that for every positive square free integer $n$, every group in $\pi^t_A(U_n)$ has a generating set with at most $\log n+C$ elements.
\end{conjecture*}
Consequences of Theorem \ref{maintheorem} in Section $1$ give some evidence for this conjecture and also generalizes Proposition $2.17$ in \cite{Ha} assuming the Galois group is solvable. In addition, this theorem also gives the following corollary, which gives a complement of a result of Jensen and Yui in \cite{JY}. They showed that any dihedral extension with Galois group $D_{2n}$ over $\Q$ ramified only at one regular prime $p$, with $p\equiv 1\,(\!\!\!\mod 4)$, has degree prime to $p$.
\begin{corc*}
Suppose $p\equiv 1\,(\!\!\!\mod 4)$ is a regular prime such that the class number of $\Q(\sqrt{p})$ is $1$. Then there are no non-abelian dihedral groups in $\pi_A(U_p)$.
\end{corc*}
For example, in the range $2\leq p\leq 100$, the primes $p=5, 13, 17, 29, 41, 53, 61, 73, 89, 97$ satisfy the conditions above. Furthermore, applying Theorem \ref{maintheorem} to the prime $p=3$ gives the following two corollaries, which are also related to results for the prime $2$ (Theorem $2.20$ and Theorem $2.23$) in \cite{Ha}. 
\begin{cord*}\label{order}
If $G$ is a solvable group in $\pi_A(U_p)$, where $p=3$, then either $G$ is cyclic, or $G/p(G)\cong \Z/2$, or $G$ has a cyclic quotient of order $27$.
\end{cord*}
\begin{core*}\label{ramificationindex}
Suppose $K/\Q$ is a Galois extension with nontrivial Galois group $G$, ramified only at the prime $p=3$ and possibly at $\infty$, with ramification index $e$. Then $9|e$ unless $G/p(G)\cong \Z/2$ or $G\cong \Z/3$.
\end{core*}
In Section $2$, we deal with the non-solvable case. We use the upper bound and lower bound of the discriminant to rule out some non-solvable groups. The main result is the following, which generalizes a result for $p=2$ in \cite{Ha}.
\begin{prop*}\label{lowerboundnonsolvable}
Let $2\leq p<23$ be a prime number. If $G\in\pa{p}$ and $|G|\leq 300$, then $G$ is solvable.
\end{prop*}
\section{Solvable extensions over $\Q$}
In this section, we will give some conditions on solvable groups that can occur as Galois groups over $\Q$ ramified only at one finite prime. 
A consequence of Harbater's Conjecture would be that if $G\in\pi_A(U_p)$ for some prime $p$, without assuming the ramification to be tame, then $G/p(G)$ is generated by at most $\log(p)+C$ elements. Thus if $p$ is very small, we expect $G$ to be very close to being a quasi-$p$ group. In fact, this holds when $p<23$ as seen the Corollary of 2.7 \cite{Ha}, i.e.\ $G/p(G)$ is cyclic of order dividing $p-1$. The following theorem is a generalization of this idea and of Proposition $2.17$ in \cite{Ha}, but with an extra assumption on solvability.
\begin{thm}\label{maintheorem}
Let $K$ be a finite Galois extension of $\Q$ ramified only at a single finite prime $p>2$, with the Galois group $G=\gal{K}{\Q}$ solvable. Let $K_0/\Q$ be an intermediate abelian extension of $K/\Q$. Let $N=\gal{K}{K_0}$ and $p(N)$ be the quasi $p$-part of $N$.
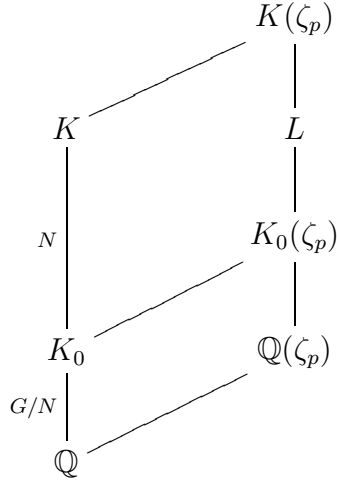
$$
\xymatrix{
 & & {K(\zeta_p)} \ar@{-}[dll] \ar@{-}[d]\\
{K} \ar@{-}[dd]_{N} & & {L}\ar@{-}[d]\\
 &  & {K_0(\zeta_p)}\ar@{-}[dll] \ar@{-}[d] \\
 {K_0} \ar@{-}[d]_{G/N}  & & {\Q(\zeta_p)} \ar@{-}[dll] \\
 {\Q} & & \\
 }
$$
Then either
\begin{enumerate}
\item[(i)] $N/p(N) \subset \Z/(p-1)$; or
\item[(ii)] there is a non-trivial abelian unramified subextension
$L/K_0(\zeta_p)$ of $K(\zeta_p)/K_0(\zeta_p)$ of degree prime to $p$ with $L$ Galois over $\Q$.
\end{enumerate}
\end{thm}
We will first give some corollaries, then prove a lemma and a proposition that will be used in the proof of Theorem \ref{maintheorem}, given at the end of this section.
\begin{remarks} \label{remarksofmaintheorem}$\quad$ 
\begin{enumerate}
\item[I)]If we let $K_0=\Q$ and $p<23$, Theorem \ref{maintheorem} is just Corollary $2.7$ \cite{Ha} in the solvable case.
\item[II)]In fact we will show in the proof of Theorem \ref{maintheorem} that the condition $(i)$ can be replaced by the condition $K/\Q$ is a quasi-$p$ extension of a totally ramified extension.
\end{enumerate}
\end{remarks}
As a direct consequence of \ref{maintheorem}, we have:
\begin{cor}
Let $K/\Q$ be a sovable Galois extension ramified only at a prime $p$ and possibly at $\infty$. Suppose $K_0=\Q(\zeta_{p^n})$ is a sub-extension of $K/\Q$ with Galois group $\gal{K}{K_0}= G$ and the class number of $K_0(\zeta_{p})$ is $1$. Then $G/p(G)$ is cyclic of order dividing $p-1$.
\end{cor}
\begin{proof}Apply Theorem \ref{maintheorem} to the subextension $K_0/\Q$. Since the class number of $K_0(\zeta_p)$ is $1$, the condition $(ii)$ in Theorem \ref{maintheorem} does not hold. So condition $(i)$ holds, i.e. $G/p(G)$ is cyclic of order dividing $p-1$.
\end{proof}

Using Theorem \ref{maintheorem}, we can prove Corollary C.
\begin{proofC*}
Suppose that $K/\Q$ is a Galois extension with group $D_{2n}$ of order $2n$, ramified only at a finite prime $p$ and possibly at $\infty$. Denote by $K_0$ the fixed field of the cyclic subgroup $\Z/n< D_{2n}$. By Theorem $1.2.2$ in \cite{JY}, we know $n$ is not divisible by $p$. Now apply Theorem \ref{maintheorem}. By the assumption the class number of $K_0$ is $1$, we know that the condition $(ii)$ in Theorem \ref{maintheorem} fails. By the second remark above we have $K/\Q$ is totally ramified, since $p$ does not divide the order of $D_{2n}$. So $\gal{K}{\Q}\cong P\rtimes C$, where $P$ is a $p$-group and $C$ is a cyclic group. We know $P$ has to be trivial, again since $p\nmid 2n$. Thus $\gal{K}{\Q}$ is cyclic.
\end{proofC*}
\begin{lem}\label{ppowermaintheorem}
Under the hypotheses of Theorem \ref{maintheorem}, if $K_0$ is a maximal $p$-power Galois sub-extension of $K/\Q$, then the condition $(i)$ can be replaced by the condition that either $G$ is a cyclic $p$-group or $N/p(N)$ is a nontrivial subgroup of $\Z/(p-1)$.
\end{lem}
\begin{proof}It suffices to show that if $N$ is a quasi-$p$ group, then $G$ is a cyclic $p$-group. So assume $G$ is not a cyclic group (in particular $G$ is non-trivial). The Galois group $\gal{K_0}{\Q}$ is cyclic, say of order $p^n$ for some $n$, since any finite $p$-group in $\pi_A(U_p)$ is cyclic (see Theorem $2.11$ in \cite{Ha}). By Class Field Theory $K_0/\Q$ is the unique cyclic sub-extension of degree $p^n$ in $\Q(\zeta_{p^{n+1}})$ since $p>2$. If $K_0'$ is another maximal $p$-power sub-extension of $K/\Q$, then $K_0'=K_0$ by the same argument. So $K_0$ is the unique maximal $p$-power sub-extension of $K/\Q$. Denote by $N$ the Galois group $\gal{K}{K_0}$. Then $N$ is normal in $G$ and it is the minimal subgroup of $G$ with index a power of $p$, since it corresponds to the unique maximal $p$-power sub-extension $K_0$. We know $N$ is nontrivial, since $G$ is not a $p$-group by assumption. Now $G$ is a nontrivial solvable group, so $G$ has a normal subgroup $\bar{N}\subset N$ such that $N/\bar{N}$ is of the form $(\Z/q)^n$ for some prime $q$ and some integer $n\geq 1$. We know $q\neq p$ from the minimality of $N$. So $N$ is not a quasi-$p$-group since every $p$-subgroup of $N$ is contained in the proper subgroup $\bar{N}$.
\end{proof}
The above lemma gives evidence for the conjecture in the introduction. Next we will apply Lemma \ref{ppowermaintheorem} to the prime $3$ to get Corollary D and consequently Corollary E.
\begin{proofD*}
Let $K/\Q$ be a solvable Galois extension ramified only at $p=3$ with Galois group $G=\gal{K}{\Q}$. Take $K_0/\Q$ to be a maximal $p$-power Galois sub-extension of $K/\Q$. The same argument as in Lemma \ref{ppowermaintheorem} shows that the Galois group $\gal{K_0}{\Q}$ is cyclic and $K_0$ is the unique maximal $p$-power Galois sub-extension of $K/\Q$. If $G$ is not cyclic, then $N$ is not quasi-$p$ by Lemma \ref{ppowermaintheorem}; hence $N/p(N)$ is non-trivial. Now suppose $G$ is not cyclic and $G/p(G)\ncong \Z/2$ and apply Theorem \ref{maintheorem}. Since the condition $(i)$ in the Theorem does not hold, the condition $(ii)$ has to hold; thus the class group of $K_0(\zeta_p)$ is nontrivial. So $K_0(\zeta_p)$ contains the cyclotomic field $\Q(\zeta_{81})$, thus $\mid\! \gal{K_0}{\Q}\!\mid \,\geq 27$, i.e. $G$ has a cyclic quotient of order $27$.
\end{proofD*}
\begin{proofE*}
If $G$ is solvable, by Corollary D we know $27\mid e$ unless $G/p(G)\cong \Z/2$ or $G$ is cyclic. In the case $G$ is cyclic, we know by class field theory $K/\Q$ is totally ramified, so $e=n=|G|$, i.e. either $e=n=3$ or $9|e$. If $G$ is non-solvable, it has order $\geq 60$. On the one hand, we know $\mid\!\mathfrak{d}_{K/\Q}\!\mid^{\frac{1}{n}} \geq 12.23$ from the discriminant table (page $400$ in \cite{Od}) for extensions of degree $\leq 60$; on the other hand, considering the discriminant upper bound (Theorem $2.6$, Chapter III, \cite{Ne}), we have $\mid\!\mathfrak{d}_{K/\Q}\!\mid^{\frac{1}{n}}\leq 3^{1+v_3(e)-\frac{1}{e}} < 3^{1+v_3(e)}$. Combining these two inequalities gives $12.23\leq \mid\!\mathfrak{d}_{K/\Q}\!\mid^{\frac{1}{n}} < 3^{1+v_3(e)},$ thus $v_3(e)\geq 2$ and $9\mid e$.
\end{proofE*}
\begin{rem*}
Corollary E does not assume the solvability of $\gal{K}{\Q}$.
\end{rem*}
For the proof of Theorem \ref{maintheorem}, we first need a lemma and a proposition.
\begin{lem}\label{grouplemma}
Let $G=P\rtimes\Z/(l_1l_2)$  be a semidirect product of a $p$-group $P$ by a cyclic group $\Z/l_1l_2$, with $p,l_2$ distinct primes and $p\nmid l_1$. Denote by $s$ the highest power of $l_2$ which divides $l_1$, i.e.\ $l_2^s\mid\mid l_1$. Suppose $G$ has a normal subgroup $N\cong \Z/l_2^{s+1}$ with the quotient group $G/N \cong \Z/(l_1l_2^{-s}p^m)$. Then $G=\Z/p^m\times \Z/(l_1l_2)$.
\end{lem}
\begin{proof}
Let $\theta:\Z/(l_1l_2)\longrightarrow \Aut(P)$ be the homomorphism corresponding to the semi-direct product $G=P\rtimes\Z/(l_1l_2)$, which sends an element $a\in\Z/(l_1l_2)$ to an automorphism $\theta_a\in \Aut(P)$. Since the $l_2$-Sylow subgroup $N$ is normal in $G$, it is the unique $l_2$-Sylow subgroup by Sylow's theorem. Identify $\Z/l_2^{s+1}$ with the subset of $G=P\rtimes \Z/(l_1l_2)\cong P\rtimes (\Z/l_1l_2^{-s}\times \Z/l_2^{s+1})$, consisting of all pairs of the form $(1,b)$ with $b\in\Z/l_2^{s+1}$. We claim $\Z/l_2^{s+1}$ acts trivially on $P$ in $G$. Now for any $(k,a)\in P\rtimes \Z/(l_1l_2)$, we have
$$\begin{tabular}{ll}
{$(k,a)(1,b)(k,a)^{-1}$} & {$=(k,ab)((\theta_{a^{-1}}(k))^{-1},a^{-1})$}\\
 & {$=(k\theta_{ab}((\theta_{a^{-1}}(k))^{-1}),b)$} \\
 & {$=(k\theta_{ba}(\theta_{a^{-1}}(k^{-1})),b)$} \\
 & {$=(k\theta_{b}(k^{-1}),b).$} \\
\end{tabular}$$
Since $\Z/l_2^{s+1}=N\vartriangleleft\ G$ by the assumption, we know $(k,a)(1,b)(k,a)^{-1}$ is of the form $(1,b)$. So $\theta_b(k^{-1})=k^{-1},\forall k\in\Z/p^m$, i.e. $\theta_b$ is trivial for all $b\in \Z/l_2^{s+1}=N$.

Next we will show the isomorphism
\begin{equation}\label{isomorphism}
P\rtimes_{\theta}(\Z/(l_1l_2^{-s})\times\Z/l_2^{s+1})\cong(P\rtimes_{\theta'}\Z/(l_1l_2^{-s}))\times\Z/l_2^{s+1} \end{equation}
where the homomorphism $\theta':\Z/l_1l_2^{-s}\longrightarrow \Aut(P)$ is the restriction of $\theta$ onto $\Z/(l_1l_2^{-s})$. On the one hand the left hand side and right hand set of \ref{isomorphism} are the same as underlying sets; on the other hand we consider the binary operation in each group. Pick any two elements $(a_1,b_1,c_1),(a_2,b_2,c_2)\in\Z/p^m\rtimes(\Z/(l_1l_2^{-s})\times\Z/l_2^{s+1})$. We have
$$\begin{tabular}{ll}
{$(a_1,b_1,c_1)(a_2,b_2,c_2)$} & {$=(a_1,(b_1,c_1))(a_2,(b_2,c_2))$}\\
 & {$=(a_1\theta_{(b_1,c_1)}(a_2),(b_1b_2,c_1c_2))$} \\
 & {$=(a_1\theta_{(b_1,c_1)}(a_2),b_1b_2,c_1c_2)$.} \\
\end{tabular}$$
And if we pick any two elements $(a_1,b_1,c_1),(a_2,b_2,c_2)\in(\Z/p^m\rtimes\Z/l_1l_2^{-s})\times\Z/l_2^{s+1}$,
$$\begin{tabular}{ll}
{$(a_1,b_1,c_1)(a_2,b_2,c_2)$} & {$=((a_1,b_1)(a_2,b_2),c_1c_2)$}\\
 & {$=((a_1\theta'_{b_1}(a_2),b_1b_2),c_1c_2)$} \\
 & {$=(a_1\theta'_{b_1}(a_2),b_1b_2,c_1c_2)$.} \\
\end{tabular}$$
Since $\Z/l_2^{s+1}$ acts trivially on $P$, we know $\theta_{(b_1,c_1)}(a_2)=\theta'_{b_1}(a_2)$. So the left hand side and right hand side of \ref{isomorphism} have the same binary operations. We conclude isomorphism \ref{isomorphism}. Now we consider the quotient group $G/(\Z/l_2^{s+1})$. By the assumption it is isomorphic to $\Z/(l_1l_2^{-s}p^m)$. So by isomorphism \ref{isomorphism} we have
$$\Z/(l_1l_2^{-s}p^m)\cong G/(\Z/l_2^{s+1})\cong P\rtimes_{\theta'}\Z/l_1l_2^{-s}.$$ So $G\cong(P\rtimes_{\theta'}\Z/(l_1l_2^{-s}))\times\Z/l_2^{s+1}\cong\Z/(l_1l_2^{-s}p^m)\times\Z/(l_2^{s+1})\cong \Z/p^m\times \Z/(l_1l_2)$.
\end{proof}
\begin{prop}\label{mainproposition}
Let $K$ be a finite solvable Galois extension of $\Q$ ramified only over one finite prime $p>2$, and let $M/\Q$ be a proper abelian subextension of $K/\Q$ such that $p\nmid |\gal{K}{M}|$. Assume there is no non-trivial abelian unramified extension of $M(\zeta_p)$ of degree prime to $p$ which is contained in $K(\zeta_p)$ and is Galois over $\Q$. Then there is a proper sub-extension $M_1/M$ in $K/M$ such that $M_1/\Q$ is abelian.
\end{prop}
\begin{proof}
Since $\gal{M}{\Q}$ is abelian and ramified only at $p$, we know by class field theory that $\gal{M}{\Q}$ is a subgroup of $\Z/p^m\times \Z/(p-1)$. Write $\gal{M}{\Q}=\Z/(p^ml_1)$ with $l_1\mid p-1$. Let $N=\gal{K}{M}$. Since $N$ is solvable, being a normal subgroup of the solvable group $G$, there is a normal subgroup $N_0$ of $N$ such that $N/N_0\cong \Z/l_2$, for some prime $l_2$ such that $(l_2,p)=1$. Let $M_0$ be the fixed field of $N_0$ in $K/M$, so $\gal{M_0}{M}\cong \Z/l_2$. Let $M_1$ be the Galois closure of $M_0$ over $\Q$, so $\gal{M_1}{M}$ is a minimal normal subgroup of a solvable group $\gal{M_1}{\Q}$. From page $85$ of \cite{Rot}, we know that $\gal{M_1}{M}\cong(\Z/l_2)^t$ for some $t\geq 1$. So the Galois group $\gal{M_1}{M_0}\cong (\Z/l_2)^{t-1}$.
$$ \xymatrix{
 K \ar@{-}[r] &  M_1=(K^{N_0})^{Gal} \ar@{-}[rr]^{(\Z/l_2)^{t-1}} & & {M_0=K^{N_0}} \ar@{-}[rr]^{\Z/l_2} & & M  \ar@{-}[rr]^{\Z/(l_1p^m)} & &  \Q} \\
$$

Pick a prime $\mathfrak{p}$ of $M_1$ over the prime $p$ of $\Q$, let $I_0\subset \gal{M_1}{M}$ be the inertia group of $\mathfrak{p}$ in $M_1$ over $M$. Since $\gal{M_1}{M}\cong (\Z/l_2)^t$ is abelian, its subgroup $I_0$ is normal and the quotient by $I_0$ is abelian. So the fixed field $M_{1,0}=M_1^{I_0}$ of $I_0$ in $M_1/M$ is unramified over $M$ at the prime $\mathfrak{p}\cap\mathcal{O}_{M_{1,0}}$, thus $M_{1,0}$ is an unramified extension of $M$ contained in $M_1$. Let $\bar{M}_{1,0}$ be the Galois closure of $M_{1,0}$ over $\Q$. So $\bar{M}_{1,0}$ is contained in $M_1$ and unramified over $M$, being the composite of unramified extensions (the conjugates of $M_{1,0}$) of $M$. And $\bar{M}_{1,0}$ is abelian over $M$ since it is contained in $M_1$. 
The extensions $\bar{M}_{1,0}/M$ and $M(\zeta_p)/M$ are disjoint since $\bar{M}_{1,0}/M$ is unramified at $p$ and $M(\zeta_p)/M$ is totally ramified at $p$. Therefore if $\bar{M}_{1,0}/M$ is a non-trivial extension, $\bar{M}_{1,0}(\zeta_p)/M(\zeta_p)$ is also non-trivial. Also we know $M(\zeta_p)=\Q(\zeta_{p^{m+1}})$ since $\gal{M}{\Q}\cong \Z/(l_1p^m)$. So $\bar{M}_{1,0}(\zeta_p)$ is a non-trivial abelian unramified extension of $\Q(\zeta_{p^{m+1}})$ of degree prime to $p$ such that $\bar{M}_{1,0}(\zeta_p)\subset K(\zeta_p)$ and $\bar{M}_{1,0}$ is Galois over $\Q$, contrary to the assumption.
$$ K \supset M_1 \supset \bar{M}_{1,0}=M_{1,0}^{Gal} \supset M_{1,0}=M_1^{I_0} \supset M  \supset \Q$$
So actually $\bar{M}_{1,0}=M$, and so $I_0=\gal{M_1}{M}\cong (\Z/l_2)^t$. But the inertia group $I_0$ is cyclic (see Corollary $4$, Page $68$ of \cite{Se}), because $M_1$ is at most tamely ramified at $\mathfrak{p}$ over $M$ as the degree of the extension $M_1/M$ is prime to $p$. So $t=1$, and the field $M_1$ is totally ramified over $M$ at $\mathfrak{p}$ with $\gal{M_1}{M}\cong \Z/l_2$. It follows $M_1$ is totally ramified over $\Q$ at the prime $p$, since $M/\Q$ is abelian thus totally ramified at $p$. So $\gal{M_1}{\Q}$ is isomorphic to the inertia group $I\cong P\rtimes C$ of $M_1$ over $\Q$ at $\mathfrak{p}$, where $P$ is a $p$-group and $C$ a cyclic group of order prime to $p$. So $C$ is a cyclic group of order $l_1l_2$, thus $I\cong \Z/p^m\rtimes \Z/l_1l_2$ with $l_1,l_2$ relatively prime to $p$. On the other hand, let $l_2^s$ be the highest power of $l_2$ which divides $l_1$, so $s\geq 0$. Consider the invariant field $M^{\Z/l_2^s}$ of $\Z/l_2^s\subset \gal{M}{\Q}$ in $M$. 
$$
\xymatrix{
K \ar@{-}[r] & M_1 \ar@{-}[rr]^{I_0=\Z/l_2} & &M \ar@{-}[rr]^{\Z/(l_2^s)} & & M^{\Z/l_2^s} \ar@{-}[rr]^{\Z/(l_1l_2^{-s}p^m)} & & \Q}
$$
It is Galois over $\Q$ since $M/\Q$ is abelian, so $\gal{M_1}{M^{\Z/l_2^s}}\vartriangleleft\gal{M_1}{\Q}$. Since $M_1/M^{\Z/l_2^s}$ is totally ramified, and tamely ramified, $\gal{M_1}{M^{\Z/l_2^s}}$ is a cyclic group. So $\gal{M_1}{M^{\Z/l_2^s}}\cong \Z/l_2^{s+1}$, and the quotient group $I/(\Z/l_2^{s+1})\cong \Z/(l_1l_2^{-s}p^m)$. It follows from the Lemma \ref{grouplemma} that $I\cong\Z/(p^m)\times \Z/l_1l_2$. So $M_1/\Q$ is abelian and $M_1\neq M$.
\end{proof}
Now we can give the proof for Theorem \ref{maintheorem}:
\begin{proofthm*}
The quasi-$p$ part $p(N)$ of $N$ is normal in $G$, since it is characteristic in the normal subgroup $N\lhd G$. Replacing $G$ and $N$ by $G/p(N)$ and $N/p(N)$ respectively, we may assume $N$ has degree prime to $p$. We will show either $K/K_0$ is cyclic of order dividing $p-1$, or $(ii)$ holds. If $K/\mathbb{Q}$ is abelian, then $K$ is inside some cyclotomic field $\Q(\zeta_{p^n})$ and $N$ is a subgroup of $\Z/p^{n-1}\times\Z/(p-1)$ and of order prime to $p$, thus $N$ is cyclic of order dividing $p-1$. Now we may assume $K_0\ne K$ and $K/\Q$ is non-abelian.

First suppose that $K/K_0$ is totally ramified. Then $K/\Q$ is totally ramified, since $K_0/\Q$ is totally ramified at $p$, so $G\cong P\rtimes C$ with $P$ a $p$-group and $C$ a subgroup of $\Z/(p-1)$. Since $K_0/\Q$ is Galois, $N\triangleleft G$. So $N\subset C$ and is cyclic of order dividing $p-1$.

Otherwise, $K/K_0$ is not totally ramified. Assume $(ii)$ doesn't hold. Let $M$ be the fixed field of $G/[G,G]$ in $K/\Q$. Since $G$ is solvable, we have $K\neq M$ and $M/\Q$ is the maximal abelian subextension in $K/\Q$. We now apply Proposition \ref{mainproposition}. So there is a subfield $M_1\neq M$ in $K/M$ such that $M_1/\Q$ is abelian, contradicting the maximality of $M$. 

We now justify Part II of Remarks \ref{remarksofmaintheorem}. After replacing $G$ and $N$ by $G/p(N)$ and $N/p(N)$ respectively, it suffices to show $K/K_0$ is totally ramified since $K_0/\Q$ is abelian thus totally ramified by class field theory. If $K/K_0$ is not totally ramified, the condition $(ii)$ holds by above.
\end{proofthm*}

\section{Nonsolvable Extensions over $\Q$}
In this section, we will prove the proposition in the introduction. We will start by considering non-abelian simple groups, which form an extreme sub-class of the non-solvable groups.  
\begin{lem}\label{pdivids}
Let $2\leq p< 23$ be a prime, and $G\in\pi_A(U_p)$ with $G$ non-abelian. Then $p||G|$; furthermore, if $G$ is simple, then $G$ is a quasi $p$-group.
\end{lem}
\begin{proof}If $p\nmid |G|$, then the quasi $p$-part $p(G)$ of $G$ is trivial since it is generated by all $p$-Sylow subgroups of $G$. By Corollary $2.7$ in \cite{Ha}, we know $G=G/p(G)$ is cyclic of order dividing $p-1$. Contradiction; thus $p\mid |G|$. Now if $G$ is simple, then $p(G)\lhd G$ implies $p(G)=G$, i.e. $G$ is a quasi $p$-group.
\end{proof}
We can use above lemmas together with the Odlyzko discriminant bound to show various simple groups cannot be in $\pi_A(U_p)$:
\begin{examples}
For $2\leq p<23$, we consider $A_5$, $S_5$, $\SL(3,2)$.
\begin{enumerate}
\item[$\bullet$]$\SL(3,2)\notin\pa{p}$ for $2\leq p<23$.
\begin{proof}
The group $\SL(3,2)$ is of order $168=2^3\cdot 3\cdot 7$. When $p\neq 2, 3,7$, if we assume $G\in \pi_A(U_p)$, by Lemma \ref{pdivids} we would have $p\mid |G|$, contradiction. In the case $p=2$, Harbater showed $\SL(3,2)\notin \pi_A(U_2)$ (Example $2.21(c)$, \cite{Ha}). In the case $p=7$, Brueggeman showed $\SL(3,2)\notin \pi_A(U_7)$ in Theorem $4.1$ \cite{Br}. For the case $p=3$, we assume $G\in \pi_A(U_3)$. Let $L/\Q$ be a corresponding Galois extension and let $e$ be the ramification index of the prime above $p$. Applying the discriminant upper bound(Theorem $2.6$, Chapter III, \cite{Ne}), we get $\dis{L}{168}\leq 3^{1+v_3(e)-1/e}$. The largest power of 3 dividing $|\SL(3,2)|=168$ is $3$, so $v_3(e)\leq 1$, thus $\dis{L}{168}\leq 3^{1+v_3(e)-1/e}\leq3^{1+1}=9.$ On the other hand, by the Odlyzko discriminant bound (Table $1$, \cite{Od}), $\dis{L}{168}\geq 15.12$ when the degree of the extension is at least $160$. Contradiction.
\end{proof}
\item[$\bullet$]The alternating group $A_5\notin\pa{p}$ for $2\leq p<23$.
\begin{proof}
The group $A_5$ is of order $60=2^2\cdot3\cdot 5$. When $p\ne 2,3,5$, by Lemma \ref{pdivids}, we know $G\in\pa{p}$ would imply $p\mid |G|$, contradiction. For $p=2$, Harbater showed that $A_5\notin\pa{p}$ (Example $2.21(a)$, \cite{Ha}). For $p=5$, we know $A_5\notin \pa{5}$ from the table \cite{Jo}. For $p=3$, we assume the simple group $A_5$ lies in $\pa{3}$ and let $L/\Q$ be a corresponding Galois extension. Applying the discriminant upper bound, we have $\dis{L}{60}\leq 3^{1+v_3(e)-1/e}$. Since the largest power of 3 dividing $|A_5|=60$ is 3, we get $v_3(e)\leq 1$, thus $\dis{L}{60}\leq 3^{1+v_3(e)-1/e}\leq3^{1+1}=9.$ On the other hand, by the Odlyzko discriminant bound (Table $1$, \cite{Od}), $\dis{L}{60}\geq 12.23$ when the degree of the extension is at least $60$. Contradiction.
\end{proof}
\item[$\bullet$]The symmetric group $S_5\notin\pa{p}$ for $2\leq p<23$.
\begin{proof}
The group $S_5$ is of order $120=2^3\cdot 3\cdot 5$. When $p\ne 2,3,5$, by Lemma \ref{pdivids} we know $G\notin \pa{p}$, for otherwise we would have $p\mid |G|$, contradiction. For $p=2$, Harbater showed $S_5\notin\pa{p}$ (Example $2.21(a)$, \cite{Ha}). For $p=5$, we know $S_5\notin \pa{p}$ from the table \cite{Jo}. For $p=3$, similarly as $A_5$, we have $\dis{L}{120}\leq 3^{1+v_3(e)-1/e}\leq3^{1+1}=9.$ But by the Odlyzko bound, we have $\dis{L}{120}\geq 14.38$ for extensions of degree at least $120$, this is a contradiction.
\end{proof}
\end{enumerate}
\end{examples}
Now we are ready to prove the Proposition in the introduction using above examples.
\begin{proofprop*}
Assume there exist non-solvable groups $G\in \pa{p}$ with order $\leq 300$, and let $G$ be such a group of smallest order. Pick a nontrivial normal subgroup $N$ of $G$. The quotient group $G/N$ is also in $\pa{p}$ but with smaller order, hence solvable. We know $N$ is non-solvable, so the order of the group $N$ is at least $60$. So $|G/N|\leq 5$, thus $G/N$ is abelian. By Lemma $2.5$ in \cite{Ha} we know $G$ is isomorphic to either $A_5$, $S_5$ or $\SL(3,2)$. By examples above, these groups do not lie in $\pa{p}$ for $2\leq p<23$.
\end{proofprop*}
\section{Tamely ramified covers of the affine line over $\F_q$}
In this section, we will denote by $k$ the rational function $\F_q(t)$, and denote by $\mathfrak{f}$ the ideal generated by an irreducible polynomial $f\in\F_q[t]$. Let $U_{\frak{f}}=\A_{\F_q}^1-(f=0)$.
\begin{prop}\label{conjugaterelation}
Let $K$ be the function field of a geometric Galois cover of the affine line over $\F_q$ with Galois
group $G$ and ramified only at a finite prime $\frak{f}$ and
possibly at $\infty$, with all ramification tame. Then there exist
$x_1,x_2,...,x_d,x_{\infty}\in G$ such that
$\langle x_1,...,x_d,x_{\infty}\rangle=G$ and $x_1...x_dx_{\infty}=1$ with
$x_1^q\sim x_2,...,x_d^q\sim x_1$ and $x_{\infty}^q\sim x_{\infty}$ (i.e.\ conjugate in $G$). Moreover, the order of each of $x_1,\dots,x_d$ is equal to the ramification index over $\mathfrak{f}$, and the order of $x_{\infty}$ is the ramification index at $\infty$. So if $K/\F_q(t)$ is unramified at $\infty$, then $x_{\infty}=1$.
\end{prop}
\begin{proof}
Suppose $\deg(\mathfrak{f})=d$. After the base change to $\F_{q^d}(t)$, the prime $\frak{f}$ splits into $d$ primes $\frak{f}_1,...,\frak{f}_d$ with degree $1$, which correspond to $d$ finite places $P_1,...,P_d$ of $\F_{q^d}(t)$. Since $\infty$ has degree $1$ in $\F_q(t)$, there is a unique place $P_{\infty}$ of $\F_{q^d}(t)$ above $\infty$.
For each place $P_i$, where $1\leq i\leq d$ or $i=\infty$, there are $g$ (independent on $i$) places $Q_{i,1},Q_{i,2},...,Q_{i,g}$ of $\bar{K}$ above $P_i$, since $\bar{K}$ is Galois over $\F_q(t)$. Each place $Q_{i,j}$ has an inertia group $I_{i,j}$. Since the extension is tamely ramified, each $I_{i,j}$ is cyclic, say generated by $x_{i,j}$, i.e. $I_{i,j}=\langle x_{i,j}\rangle$. Fixing $i$, the inertia groups $I_{i,j}$ are all conjugate in $G$. The Galois group $\gal{F_{q^d}(t)}{F_q(t)}\cong\Z/d=\langle\sigma\rangle$ is generated by the Frobenius map $\sigma$, which cyclicly permutes the places $P_{i}$ where $1\leq i\leq d$; say $\sigma(P_i)=P_{i+1}$ with $(i$ mod $d)$. Also, $\sigma(P_{\infty})=P_{\infty}$. On the other hand, there is a choice of places $Q_i$ above $P_i$ for $1\leq i\leq d$ and $i=\infty$ such that the generators $x_i$ of the corresponding inertia groups generate the Galois group $G$, i.e.\ $\langle x_1,...,x_d,x_{\infty}\rangle =G$, and $x_1x_2...x_dx_{\infty}=1$. (Namely, these $Q_i$'s are specializations of corresponding ramification points of a lift of this tame cover to characteristic $0$, as in \cite{Gro} XIII.) Since all the places $P_i$ with $1\leq i\leq d$ lie over the same closed point $\frak{f}$, there is an additional condition on the group $G$. Namely, the Frobenius map $\sigma$ takes the place $Q_1$ to some $Q_2'$ over $P_2$; but the inertia group $I_2$ and $I_2'$ are conjugate, so $x_1^q\sim x_2$. Similarly $x_2^q\sim x_3,...,x_d^q\sim x_1$. Since the Frobenius map $\sigma$ maps the place $Q_{\infty}$ to some place $Q'_{\infty}$ over $P_{\infty}$, we have $x_{\infty}^q\sim x_{\infty}$.
\end{proof}
Now we consider Galois covers of the projective line $\PP^1_{\F_q(t)}$ ramified only at a finite prime $\mathfrak{f}$, generated by an irreducible polynomial $f$ in $\F_q[t]$, and unramified at $\infty$. So $U_{\mathfrak{f}}=\PP^1_{\F_q(t)}-(f=0)$.
\begin{cor}\label{semidirectproduct}
Let $n$ be a positive integer, and let $q$ be a power of an odd prime or $q=2$ or $4$. Suppose the degree $d$ of the prime $\mathfrak{f}$ is not divisible by $n$, and consider the group $G\cong P\rtimes \Z/(q^n-1)$, where $P$ is a $p$-group of order $q^n$ and the semi-direct product $G$ corresponds to the action of $\F^*_{q^n}$ on $\F_{q^n}$. Then $G\notin\pi_A^t(U_{\mathfrak{f}})$ where $U_{\mathfrak{f}}=\PP^1_{\F_q(t)}-(f=0)$.
\end{cor}
\begin{proof}Otherwise suppose $G\in\pi_A^t(U_{\mathfrak{f}})$, where $d=\deg(f)$ is not divisible by $n$. By Proposition \ref{conjugaterelation}, we have $G=\langle x_1,\dots, x_d\rangle$ with relations $x_1^q\sim x_2,\dots,x_d^q\sim x_1$ and $x_1\cdot\cdot\cdot x_d=1$. From the relation $x_1^q\sim x_2,\dots,x_d^q\sim x_1$, we have $x_i^{q^d}\sim x_i$ for $1\leq i\leq d$. Write $x_i=(a_i, b_i)$, where $a_i\in P$ and $b_i\in \Z/(q^n-1)$. Then $x_i\notin P$, for otherwise all $x_i$'s are in the normal subgroup $P$ and can not generate the group $G$. Also one of the $b_i$ has to be a generator of $\Z/(q^n-1)$, otherwise the $x_i$'s will not generate the whole group $G$. Now $x_i^{q^d}\sim x_i$ implies $b_i^{q^d}=b_i$, i.e. $q^d\equiv 1 \,\,(\!\!\mod q^n-1)$. Since $q^d\equiv 1 \,\,(\!\!\mod q^n-1)$ if and only if $n\mid d$, we have $n|d$. This is a contradiction.
\end{proof}
\begin{example*}Let $p=2$ and $n=2$ in Corollary \ref{semidirectproduct}. We have $A_4\notin \pi_A^t(U_{\mathfrak{f}})$ for any prime generated by an irreducible polynomial $f\in\F_2[t]$ of odd degree.
\end{example*}
\begin{rem*}In fact, the conclusion in \ref{semidirectproduct} can also be obtained using cyclotomic function fields, i.e.\ there are no cyclic extensions of order $q^n-1$ over $\F_q(t)$ ramified only at one finite prime $\mathfrak{f}$  (and the ramification is tame) of degree not divisible by $n$. In fact, such extension would be inside a cyclotomic function field of degree $q^d-1$ (see Theorem $2.3$, \cite{Hay}). 
\end{rem*}
Applying the proposition to dihedral groups $D_{2n}=\langle r,s\mid r^n=s^2=1, rs=sr^{-1}\rangle$ and symmetric groups $S_n$, we get Corollary A and Corollary B as follow.
\begin{proofA*}Suppose that $K/\F_q(t)$ is a geometric Galois extension with group $D_{2n}$, ramified only at a finite prime $f$ with $\deg f=d$. Applying Proposition \ref{conjugaterelation}, we know $G=\langle x_1,\cdots,x_d\rangle $ with relations $x_1\cdots x_d=1$ and
$x_1^q\sim x_2,\cdots,x_d^q\sim x_1$. Now we divide the situation into two cases ($n$ is even and $n$ is odd):
\begin{case}[$n=2k$]The conjugacy classes in $D_{2n}$ are
$\{1\}$, $\{r^k\}$, $\{r^{\pm1}\}$, $\{r^{\pm2}\},\cdots,\{r^{\pm (k-1)}\}$, $\quad \{sr^{2b}\mid b=1,\cdots,k\}$, $ \{sr^{2b-1}\mid b=1,\cdots,k\}.$ If one of $x_1,\cdots,x_d$ is a power of $r$, then all the $x_i$'s have to be powers of $r$ because of the conjugacy relations. This is a contradiction since such $x_i$'s cannot generate the group $D_{2n}$. So we have $x_i=sr^{t_i},\,\,1\leq i\leq d$, for some integers $t_i$. If some $t_i$ is even, then again by the conjugacy relations all $t_i$'s have to be even. This is a contradiction since such $x_i$'s cannot generate the group $D_{2n}$. So we can write $x_i=sr^{2k_i+1},\,\, 1\leq i\leq d$, which implies
$sr^{2k_1+1}\cdot sr^{2k_2+1}\cdot\cdot\cdot sr^{2k_{d-1}+1}\cdot sr^{2k_d+1}=1.$
This is impossible. Therefore $D_{4k}\notin \pi_A^t{(U_f)}$.
\end{case}
\begin{case}[$n=2k+1$]The conjugacy classes in $D_{2n}$ are $\{1\}$, $\{r^{\pm1}\}$, $\{r^{\pm2}\},...,\{r^{\pm k}\},\{sr^b\mid b=1,...,n\}.$ Similarly we can write $x_i=sr^{k_i},\,\,1\leq i\leq d$, which gives $sr^{k_1}\cdot sr^{k_2}\cdot\cdot\cdot sr^{k_{d-1}}\cdot sr^{k_d}=1.$ This is impossible if $2\nmid d$. Thus $D_{4k+2}\notin \pi_A^t{(U_f)}$ if the degree of the prime $f$ is odd.
\end{case}
\end{proofA*}
\begin{proofB*}
Suppose that $K/\F_q(t)$ is a geometric Galois extension with group $S_n$, ramified only at a finite prime $f$ with $\deg{f}=d$. Applying Proposition \ref{conjugaterelation}, we know there exist $x_1,...,x_d$ such that $G=\langle x_1,...,x_d\rangle$ with relations $x_1\cdot\cdot\cdot x_d=1$ and $\hspace{2em}x_1^q\sim x_2,...,x_d^q\sim x_1.$ If $2|q$, all $x_i$'s are even permutations since two permutations are conjugate in $S_n$ if and only if they have the same cycle structure. This is impossible since they cannot generate $S_n$; If $2\nmid q$, all $x_i$'s are of the same parity. Since they generate $S_n$, they have to be odd permutations. So if $d$ is odd, the product $\prod_{i=1}^d{x_i}$ of an odd number of odd permutations $x_i$'s is still an odd permutation, which cannot be $1$, contradiction. 
\end{proofB*}
\bigskip

\providecommand{\bysame}{\leavevmode\hbox
 to3em{\hrulefill}\thinspace}

 \end{document}